\documentclass[11pt]{article}
\usepackage{amssymb}
\usepackage{latexsym}
\newtheorem{Theorem}{Theorem}[part]

\newtheorem{Proposition}{Proposition}[part]

\newtheorem{Lemma}{Lemma}[part]
\newtheorem{Corollary}{Corollary}[part]
\newtheorem{Remark}{Remark}[part]
\newtheorem{Example}{Example}[part]

\makeatletter \@addtoreset{equation}{section}

\@addtoreset{Definition}{section}

\@addtoreset{Theorem}{section}

\@addtoreset{Proposition}{section}

\@addtoreset{Property}{section}

\@addtoreset{Assumption}{section}

\@addtoreset{Corollary}{section}

\@addtoreset{Lemma}{section}

\@addtoreset{Remark}{section}

\@addtoreset{Example}{section}

\def \Max{\displaystyle\max}

\def \proof{{\noindent \bf Proof. }}
\def \ep{\hbox{ }\hfill$\Box$}
\def\reff#1{{\rm(\ref{#1})}}

\def\Ac{{\cal A}}

\def\eps{\varepsilon}

\def\Cc{{\cal C}}

\def\Dc{{\cal D}}

\def\Fc{{\cal F}}
\def\Gc{{\cal G}}

\def\Kc{{\cal K}}

\def\Lc{{\cal L}}

\def\Pc{{\cal P}}
\def\Qc{{\cal Q}}
\def\Xc{{\cal R}}

\def\Tc{{\cal T}}
\def\Xc{{\cal X}}
\def\Zc{{\cal Z}}

\def\Zc{{\cal Z}}

\def\eps{\varepsilon}

\def\no{\noindent}
\def\Pas{\mathbb{P}-\mbox{a.s.}}

\def\x{\times}

\def\05{\frac{1}{2}}
\def\-1{^{-1}}
\def\1{{1\hspace{-1mm}{\rm I}}}
\def\={\;=\;}
\def\.{\;.}

\title{  }
\author{ }

\def\be{\begin{eqnarray}}
\def\ee{\end{eqnarray}}

\def\b*{\begin{eqnarray*}}
\def\e*{\end{eqnarray*}}

\def\And{\;\mbox{ and }\;}

\def\Esp#1{\mathbb{E}\left[#1\right]}
\def\Espt#1{\mathbb{E}^{\tilde \P}\left[#1\right]}
\def\Pro#1{\mathbb{P}\left[{#1}\right]}

\def\And{\mbox{ and } }
\def\pourtout{\mbox{ for all } }

\def \E{\mathbb{E}}
\def \F{\mathbb{F}}
\def \M{\mathbb{M}}
\def \N{\mathbb{N}}
\def \R{\mathbb{R}}

\def\P{\mathbb{P}}
\def\Q{\mathbb{Q}}

\def\T{\mathbb{T}}

\def\vart{\vartheta}
\def\ubf{{\bf 1}}

\title{On the  Hedging of American Options  in Discrete Time
Markets with Proportional Transaction Costs}

\begin{document}
\vspace{-15mm}
\author{Bruno BOUCHARD\footnote{The  authors would like to
thank Y. Kabanov for fruitful discussions on the subject.}
             \\\small  Laboratoire de Probabilit{\'e}s et
             \\\small  Mod{\`e}les Al{\'e}atoires
             \\\small  CNRS, UMR 7599
             \\\small  Universit{\'e} Paris 6 and CREST
             \\\small  e-mail: bouchard@ccr.jussieu.fr
             \and
         Emmanuel TEMAM
             \\\small  Laboratoire de Probabilit{\'e}s et
             \\\small  Mod{\`e}les Al{\'e}atoires
             \\\small  CNRS, UMR 7599
             \\\small  Universit{\'e} Paris 7
             \\\small  e-mail: temam@math.jussieu.fr
             }

\date{This version : January 2005}

\maketitle

\vspace{-5mm}

\begin{abstract}
In this note, we consider a general discrete time financial market
with proportional transaction costs as in Kabanov and Stricker
\cite{KabStri}, Kabanov et al. \cite{KSRstrict}, Kabanov et al.
\cite{KSR} and Schachermayer \cite{schach}. We provide a dual
formulation for the set of initial endowments which allow to
super-hedge some American claim. We show that this extends the result  of
Chalasani and Jha \cite{CJ} which was obtained in a model with
constant transaction costs and risky assets which evolve on a
finite dimensional tree. We also provide fairly general conditions
under which the expected formulation in terms of stopping times
does not work.
\end{abstract}

\vspace{7mm}

\noindent {\bf Key words~:} Sum of random convex cones,
Transaction costs, American option.

\noindent {\bf MSC Classification (2000):} 91B28, 60G40.

\bigskip

\section{Introduction}

We consider a discrete time financial market with proportional
transaction costs. These markets have already been widely studied.
In particular, a proof of the fundamental theorem of asset pricing
was given in Kabanov and Stricker \cite{KabStri} in the case of
finite $\Omega$ and further developed in Kabanov et al.
\cite{KSRstrict}, \cite{KSR}, R$\acute{{\rm a}}$sonyi
\cite{rasThese}, Schachermayer \cite{schach} among others. In
these papers, a super-replication theorem is also provided for
European contingent claims. The aim of this paper is to extend this theorem
to American options. It is well known that, for frictionless
markets, the super-replication price of an American claim admits a
dual representation in terms of  stopping  times. However, it is
proved in Chalasani and Jha \cite{CJ} that, in markets with fixed
proportional transaction costs and with assets evolving on a
finite dimensional tree, this formulation does not hold anymore.
In their setting, they show that a correct dual formulation can
be obtained if we replace stopping times by {\it randomized}
stopping times, which amounts to work with what they call
``approximate martingale node-measures''.
 In this paper, we provide a new dual formulation
for the price of American option which works in  the general
framework of $\Cc$-valued processes as introduced in Kabanov et
al. \cite{KSRstrict}, and extends the dual formulation of
Chalasani and Jha \cite{CJ}.

\vspace{0.3cm} The rest of the paper is organized as follows. In
Section 2, we describe the model and give the dual formulation.
The link between our formulation and the one obtained by Chalasani
and Jha \cite{CJ} is explained in Section 3. Section 4 is devoted
to counter-examples.  The proof of the dual formulation is
provided in Section 5.

\section{Model and main result}

\subsection{Problem formulation}

\no Set $\T$ $=$ $\{0,\ldots,T\}$ for some $T\in \N\setminus\{0\}$
and let $(\Omega,\Fc,\P)$ be a complete probability space endowed
with a filtration $\F$ $=$ $(\Fc_t)_{t\in \T}$. In all this paper,
inequalities involving random variables have to be understood in
the $\Pas$ sens. We assume that $\Fc_T$ $=$ $\Fc$ and that $\Fc_0$
is trivial. Given an integer $d\ge 1$, we denote by  $\Kc$  the
set of $\Cc$-valued processes $K$ such that
$\R^{d}_+\setminus\{0\} \subset \mbox{int}(K_t)$  for all $t\in
\T$.\footnote{ Here, we follow Kabanov et al. \cite{KSR} and say
that a sequence of set-valued mappings $(K_t)_{t\in \T}$ is a
$\Cc$-valued process if there is a countable sequence of
$\R^{d}$-valued $\F$-adapted processes $X^n$ $=$ $(X^n_t)_{t\in
\T}$ such that, for every $t \in \T$, $\Pas$ only a finite but
non-zero number of $X^n_t$ is different from zero and $K_t$ $=$
cone$\{X^n_t\;,\;n\in \N\}$. This means that $K_t$ is the
polyhedral cone generated by the $\Pas$ finite set
$\{X^n_t\;,\;n\in \N\;\And\;X^n_t\ne 0\}$.} \vspace{0.4cm}

\no Following the modelization of Kabanov et al. \cite{KSR}, for a
given $K$ $\in \Kc$ and $x$ $\in$ $\R^d$, we define the process
$V^{x,\xi}$ by
    \b*
    V^{x,\xi}_t&:=&x+\sum_{s=0}^t \xi_s \;\;,\;t\in \T\;,
    \e*
where  $\xi$  belongs to
    \b*
    \Ac(K) =\{\;\xi \in L^0(\R^{d};\F)\;\; \mbox{s.t.}\;\; \xi_t\;\in\; -K_t&&
    \;\;\pourtout\; t\in \T\}\;,
    \e*
and, for a random set $E$ $\subset \R^d$ $\Pas$ and $\Gc\subset
\Fc$, $L^0(E;\F)$ (resp. $L^0(E;\Gc)$) is the collection of
$\F$-adapted processes (resp. $\Gc$-measurable variables) with
values in $E$.\medskip

The financial interpretation is the following: $x$ is the initial
endowment in number of physical units of some financial assets,
$\xi_t$ is the amount of physical units of assets which is
exchanged at $t$ and $-K_t$ is the set of affordable exchanges
given the relative prices of the assets and the level of
transaction costs.\medskip

Before to go on, we illustrate this modelization through an
example (see also Section 3,  Kabanov and Stricker \cite{KabStri}
and Kabanov et al. \cite{KSR}).

\begin{Example}\label{ex marche fi}  Let us a consider a currency market
with $d$ assets whose price process is modelled by the
$(0,\infty)^d$-valued $\F$-adapted process $S$. Let  $\M^d_+$
denote the set of $d$-dimensional square matrices with
non-negative entries. Let  $\lambda$ be a $\M^d_+$-valued
$\F$-adapted process and consider the $\Cc$-valued process
$(K_t)_{t\in \T}$ defined by
    \[
    \hspace{-1mm} K_t(\omega)=\left\{x\in \R^d~:~\exists a \in \M^d_+
    \mbox{ s.t. }
    x^i + \sum_{j\ne i\le d} a^{ji} - a^{ij}\pi^{ij}_t(\omega)
    \ge 0 \;\;\forall\;i\le d\right\}\;,
    \]
where $\pi^{ij}_t:=(S^j_t/S^i_t)(1+\lambda^{ij}_t)$  for all
$i,j\le d$ and $t\in \T$. In the above formulation, $a^{ij}$
stands for the number of units of asset $j$ obtained by selling
$a^{ij}\pi^{ij}_t$ units of assets $i$. $\lambda^{ij}_t$ is the
coefficient of proportional transaction costs paid in units of
asset $i$ for a transfer from asset $i$ to asset $j$. If $\xi_t
\in -K_t$, then we can find some financial transfers
$\eta_t=(\eta^{ij}_t)_{i,j\le d} \in L^0(\M^d_+;\Fc_t)$ such that
    \b*
    \xi^i_t \le \sum_{j\ne i\le d} \eta^{ji}_t -
    \eta^{ij}_t\frac{S^j_t}{S^i_t} (1+\lambda^{ij}_t )\;,
    \;\;i\le d\;\;,
    \e*
i.e. the global change in the portfolio position is the result of
single exchanges, $\eta^{ij}_t$, between the different financial
assets, after possibly throwing away some units of these assets.

The random set $K_t$ denotes the so-called {\it solvency region},
i.e. $V_t \in K_t$ means that, up to an immediate transfer
$\xi_t$, $V_t$ can be transformed into a portfolio with no {\it
short-position} $\tilde V_t=V_t+\xi_t$ $\in  \R^d_+ $.
Observe that we can assume, without loss of generality, that
    \b*
    (1+\lambda^{ik}_t)(1+\lambda^{kj}_t) &\ge& (1+\lambda^{ij}_t)
    \;\;\;i, j, k\le d,\; t\in \T\;.
    \e*
Indeed, if this condition is not satisfied then any ``optimal''
strategy would induce an effective transaction cost equal to
$\tilde \lambda^{ij}_t:=(1+\lambda^{ik}_t)(1+\lambda^{kj}_t)-1$.
\end{Example}

\medskip

\no The set of all portfolio processes with initial endowment $x$
is given by
    \b*
    A(x;K)&:=&\left\{V^{x,\xi},\; \xi \in \Ac(K)\right\}
    \e*
so that
    \b*
    A_t(x;K)&:=&\left\{V_t,\; V \in A(x;K)\right\}\;
    \e*
corresponds to the collection of their  values at time $t \in
\T$.\medskip

\no  It is known from the work of Kabanov and Stricker
\cite{KabStri}, Kabanov et al. \cite{KSR}, Kabanov et al.
\cite{KSRstrict} and Schachermayer \cite{schach}, see also the
references therein, that, under mild no-arbitrage assumptions (see $NA^s(K)$ and $NA^r(K)$ below), the set
$A_T(x;K)$ can be written as
    \begin{equation}\label{eq formu euro}
    \left\{g \in L^0(\R^d;\Fc)~:~\Esp{Z_T\cdot g - Z_0 \cdot x} \le 0
    ,\; \pourtout\;Z \in \Zc(K),\;(Z\cdot g)^- \in L^1(\R;\P) \right\}
    \end{equation}
where $\Zc(K)$ is the set of $(\F,\P)$-martingales $Z$ such that
    \b*
    Z_t \;\in\; K^*_t && \;\pourtout \;t\in \T\;,
    \e*
and $K^*_t(\omega)$ denotes the positive polar of $K_t(\omega)$,
i.e.
    \b*
    K^*_t(\omega)&:=&\left\{y\in \R^d~:~x\cdot y \ge
    0\;,\;\pourtout x\in K_t(\omega)\right\}\;.
    \e*
The operator ``$\cdot$'' denotes the natural scalar product on
$\R^d$ and $L^1(E;\P)$ (resp. $L^1(E;\F,\P)$)  is the subset of
$\P$-integrable elements of $L^0(E;\Fc)$ (resp. $L^0(E;\F)$).

\vspace{0.4cm}

\no In this paper, we are interested in
    \b*
    A^s(x;K)&:=&\left\{ \vartheta \in
    L^0(\R^d;\F)~:~  V-\vartheta \in -\Ac(K)
     \mbox{ for some } V \in A(x;K) \right\}\;,
    \e*
the set of processes which are dominated  by a portfolio in the
sense of $K$: $V_t-\vart_t \in K_t$, for all $t\in \T$.
The relation  $V_t-\vart_t \in K_t$ means that there is an
immediate financial transaction $\xi_t \in -K_t$ such that
$V_t+\xi_t=\vart_t$. Hence, $A^s(x;K)$ can be interpreted as the
set of American claims $\vartheta$, labelled in physical units of the financial assets,
 which are super-hedgeable when
starting with an initial wealth equal to $x$.\medskip

\no More precisely, our aim is to provide a dual formulation for
    \begin{equation}
    \label{GammaDefinition}
    \Gamma(\vartheta;K):=\left\{x\in \R^d~:~\vartheta \in
    A^s(x;K)\right\}\;,
    \end{equation}
the set of initial holdings $x$ that allow to super-hedge
$\vartheta$.

\subsection{Dual formulation}

\no In analogy with the standard result for markets without
transaction cost, one could expect that $\Gamma(\vartheta;K)$ can
admit the dual formulation
    \begin{equation}\label{eq mauvaise formualtion duale}
     \Theta(\vart;K)=\left\{x \in \R^d~:~\sup_{\tau \in \Tc(\T)}
     \Esp{Z_\tau \cdot \vart_\tau - Z_0 \cdot x} \le 0
    \;,\; \pourtout\;Z \in \Zc(K)\right\}
    \end{equation}
where $\Tc(\T)$ is the set of all $\F$-stopping times with values
in $\T$. If we assume that $\vartheta$ is bounded from below,
component by component, this formulation is obtained by taking the
supremum over all stopping times in \reff{eq formu euro} as in the
frictionless case. However, this characterization does not hold
true in general, as shown in Section 4.  This phenomenon was
already pointed out in Chalasani and Jha \cite{CJ} in a model
consisting of one bank account and one risky asset evolving on a
finite dimensional tree. In Chalasani and Jha \cite{CJ}, the
authors show that a correct dual formulation can be obtained if we
replace stopping times by {\it randomized stopping
times}.\vspace{0.5cm}

\no In our general framework, this amounts to introduce a new set
of dual variables, see Section 3 for an interpretation in terms of {\it randomized stopping
times}. For $\tilde \P \sim \P$, the associated set of
dual variables, $\Dc(K,\tilde \P )$, is defined as the collection
of processes $Z$ $\in$ $L^1(\R^d_+;\F,\tilde \P)$ such that
    \b*
    Z_t \in K^*_t \; \mbox{ and }\; \bar Z_t&:=& \Espt{\sum_{s=t}^T
    Z_s~|~\Fc_t} \in K^*_t \;\;\pourtout\;t\in \T\;.
    \e*

\begin{Example} In the
model of Example \ref{ex marche fi}, we have
    \b*
    K^*_t(\omega)&=&\left\{y\in \R^d_+~:~ y^jS^i_t(\omega) \le y^i S^j_t(\omega)(1+\lambda^{ij}_t(\omega))\;,  \;\;i\ne j\le
    d\right\}\;.
    \e*
It follows that  $\Dc(K,\tilde \P )$ is    the collection of
processes $Z$ $\in$ $L^1(\R^d_+;\F,\tilde \P)$ such that
    \b*
   Z^j_tS^i_t \le Z^i_t S^j_t(1+\lambda^{ij}_t) \; \mbox{ and }\; \bar Z^j_tS^i_t \le \bar Z^i_t S^j_t(1+\lambda^{ij}_t)
    \;\;\forall\;i, j\le d\;,\;t\in \T\;.
    \e*
\end{Example}

\medskip

In the following, we shall say that a subset $B$ of $L^0(\R^d;\F)$
is {\sl closed in measure} if it is closed in probability when
identified as a subset of $L^0(\R^{d \times (T+1)};\Fc)$, i.e.
    \b*
    v^n \in B \And \forall \eps >0 \quad
    \lim_{n\to \infty}\Pro{\sum_{t \in \T} \|
    v^n_t - v_t \| > \eps} = 0 &\Longrightarrow & v \in B\;.
    \e*

\no We then have the following characterization of $A^s(K)$ $:=$
$A^s(0;K)$.

\begin{Theorem} \label{thm main}
Assume that $A^s(K)$ is {\sl closed in measure} and that the
no-arbitrage condition
    \b*
    NA(K)~:~  \;\;A_T(0;K) \cap L^0(\R^{d}_+;\Fc) &=& \{0\}
    \e*
holds.  Then, for all $\tilde \P\sim \P$, there is $Z\in
\Dc(K,\tilde \P)$ with values in $(0,\infty)^d$. Moreover, the
following assertions are equivalent~:

\begin{enumerate}

\item[(i)] $\vart \in A^s(K)$

\item[(ii)] for all $\tilde \P \sim \P$ and $Z \in \Dc(K,\tilde
\P)$ such that $(\vart\cdot Z)^- \in L^1(\R;\F,\tilde \P)$ we have
    \b*
    \Espt{\sum_{t=0}^T \vart_t \cdot Z_t} &\le& 0\;.
    \e*

\item[(iii)] for some $\tilde \P \sim \P$ we have
    \b*
    \Espt{\sum_{t=0}^T \vart_t \cdot Z_t} &\le& 0\;
    \e*
for all $Z \in \Dc(K,\tilde \P)$ such that $(\vart\cdot Z)^- \in
L^1(\R;\F,\tilde \P)$.
\end{enumerate}

\end{Theorem}

\no Since $A^s(K)$ $=$ $A^s(0;K)$ $=$ $A^s(x;K)-x$, this
immediately provides a dual formulation for $\Gamma(\vart;K)$.

\begin{Corollary}
\label{cor form duale} Let the conditions of Theorem \ref{thm
main} hold. Then, for all $\vart \in L^0(\R^d;\F)$,

\small{  \[
    \Gamma(\vart;K)
    =
    \left\{ x\in \R^d~:~\Esp{\sum_{t=0}^T \vart_t \cdot Z_t} \le \bar
    Z_0 \cdot x \;\;\forall\; Z \in \Dc(K;\P),\;(Z\cdot \vart)^- \in L^1(\R;\F,\P)\right\}.
    \]}
\end{Corollary}

\begin{Remark}{\rm The integrability condition  on $(Z\cdot
\vart)^-$ is trivially satisfied if there is some $\R^d$-valued
constant $c$ such that  $\vart_t + c \in K_t$  for all $t\in
\T$, i.e. the {\sl liquidation value} of $\vartheta$ is uniformly
bounded from below. Indeed, in that case $Z_t\cdot(\vart_t + c)\ge
0$  for all $Z_t\in L^0(K^*_t;\Fc)$.}
\end{Remark}

\no Following the approach of Kabanov et al. \cite{KSRstrict} and
Kabanov et al. \cite{KSR} the closure property of  $A^s(0;K)$ can
be obtained under the general assumption
    \be
    \label{eq_ha}\;\;\xi \in \Ac(K) \And \sum_{t\in \T} \xi_t=0
     &\Longrightarrow& \xi_t\in K^0_t\;\;\;\pourtout t\in
    \T\;
    \ee
where $ {K}^0$ $=$ $( { K}^0_t)_{t\in \T}$ is defined by $ {
K}^0_t= {K}_t\cap(- {K}_t)$ for $t\in \T$.

\begin{Proposition}\label{lem closure} Assume that \reff{eq_ha} holds, then $A^s(K)$ is {\sl closed
  in measure}.
\end{Proposition}

\begin{Remark}

{\rm 1. In the case of efficient frictions, i.e. $K^0_t=\{0\}, \;
$ $ \forall t \in \T$, it is shown in   Kabanov et al.
\cite{KSRstrict} that the assumption \reff{eq_ha} is a consequence
of the {\sl strict no-arbitrage} property
    \b*
    NA^s(K):\;  A_t(0;K)\cap L^0(K_t;\Fc_t)&\subset& L^0(K^0_t;\Fc_t)
    \;\pourtout\; t\in \T \;.
    \e*
The financial interpretation of the assumption $K^0_t=\{0\}$ is that there is no
couple of assets which can be exchanged freely, without paying
transaction costs. In the model of Example \ref{ex marche fi}, it
is easily checked that it is equivalent to
$\lambda^{ij}_t+\lambda^{ji}_t>0$  for all $t\in \T$.

\vspace{3mm}

\no 2. In the case where $K^0_t$ may not be trivial,  \reff{eq_ha}
holds under the {\sl robust no-arbitrage} condition introduced by
Schachermayer \cite{schach} and further studied by Kabanov et al.
\cite{KSR},
$$
NA^r(K ):\;  NA(\tilde K ) \mbox{ holds for some  $\tilde K \in
\Kc$ which dominates $ K $,}
$$
where $\tilde K $   dominates $ K $ if $K_t\setminus   K
^0_t\;\subset\;\mbox{ri}(  \tilde K_t)    \;
\;\;\pourtout\; t\in \T\;$.

In finance, this means that there is a model with strictly bigger
transaction costs (in the directions where they are not equal to
zero) in which there is still no-arbitrage in the sense of $NA$.

\vspace{3mm}

\no 3. It is shown in Penner \cite{Penner} that the condition
$K^0_t=\{0\}$ in 1. can be replaced by the weaker one:
$L^0(K^0_t;\Fc_{t-1})\subset L^0(K^0_{t-1};\Fc_{t-1})$ for all
$1\le t\le T$. See also R$\acute{{\rm a}}$sonyi \cite{rasThese}.
 }\end{Remark}

\section{Interpretation in terms of ``approximate martingale node-measures'' and ``randomized stopping times''}

In this section, we show that the dual formulation of Corollary
\ref{cor form duale} is indeed an extension of the result of
Chalasani and  Jha \cite{CJ}. To this purpose, we consider the
example treated in the above paper. It corresponds to a financial
market with one non-risky asset $S^1$ and one risky asset $S^2$. For sake of simplicity, we restrict ourselves to this $2$-dimensional case, although it
should be clear that the above arguments can be easily extended to the multivariate setting.
We also assume  that the interest rate
associated to the non-risky asset is equal to zero and normalize
$S^1\equiv 1$, otherwise all quantities have to be divided by
$S^1$ which amounts to considering $S^1$ as a num\'eraire and
working with discounted values. Here, $S^2$ is a
$(0,\infty)$-valued $\F$-adapted process. Given $\mu$ and
$\lambda$ in $(0,1)$, the model considered in Chalasani and  Jha
\cite{CJ} corresponds to the $\Cc$-valued process $(K_t)_{t\in \T}$ defined by
    \[
    K_t(\omega)= \left\{x\in \R^2~:~
    x^1 + S^2_t(\omega)\left[(1-\mu) x^2\1_{x^2\ge 0} + (1+\lambda) x^2\1_{x^2<0} \right]
    \ge 0 \right\}\;,
    \]
so that
        \[
    K^*_t(\omega)=\left\{y\in \R^2_+~:~
         S^2_t(\omega)(1-\mu)y^1\le y^2 \le    S^2_t(\omega) (1+\lambda)y^1
    \right\}\;.
    \]
It follows that  $\Dc(K,  \P )$ is the collection of
processes $Z$ $\in$ $L^1(\R^2_+;\F, \P)$ such that
   \be\label{eq def Dc dan Cha Jha}
    S^2_t (1-\mu)Z^1_t\le Z^2_t \le    S^2_t (1+\lambda)Z^1_t \;
    \mbox{ and }\;
    S^2_t (1-\mu)\bar Z^1_t\le \bar Z^2_t \le    S^2_t  (1+\lambda)\bar Z^1_t
    \;\;
    \ee
for all $t \in \T$.

\subsection{ ``Approximate martingale node-measures''}

We first provide an alternative dual formulation in terms of  what
Chalasani and  Jha \cite{CJ} call ``approximate martingale
node-measures''. Although we are not considering a finite
probability space, we keep the term ``node measure'' used in the
above paper for ease of comparison.\medskip

Given $Z$ in $\Dc(K,  \P )$, let us define
$\hat Z^1$ $=$ $Z^1/S^1$ $=$ $Z^1$ and $\hat Z^2$ $=$ $Z^2/S^2$ so
that \reff{eq def Dc dan Cha Jha} can be written equivalently in
    \begin{equation}\label{eq def Dc dan Cha Jha equiv}
    (1-\mu)\hat Z^1_t\le \hat Z^2_t \le     (1+\lambda) \hat Z^1_t \;
    \mbox{ and }\;
     S^2_t (1-\mu) E^{q^Z}_t[{\bf 1}]  \le E^{q^Z}_t[\chi^Z S^2] \le   E^{q^Z}_t[{\bf 1}]   S^2_t   (1+\lambda)
    \;\;
    \end{equation}
where ${\bf 1}$ denotes the constant process equal to $1$, and,  for a  bounded from below  process $\alpha$
    \b*
    E^{q^Z}_t[\alpha] :=  \Esp{\sum_{t\le s\le T}q^Z_s \alpha_s~|~\Fc_t}
        &\mbox{with}&
    \chi^Z_t\;:=\; \hat Z^2_t/\hat Z^1_t  \And  q^Z_t\;:=\;\hat Z^1_t / \Esp{\sum_{s\in \T}   \hat Z^1_s }\;,
    \e*
for all $t\in \T$. Here, we use the convention $0/0=0$. Hence, any
element $Z$ in $\Dc(K,  \P )$ with $Z^1\ne 0$ is such that
$(\chi^Z,q^Z)$ belongs to the set $\Qc(K,\P)$ of elements
$(\chi,q)$ of $L^1( \R^2_+;\F, \P)$ such that, for all $t\in \T$,
     \begin{equation}\label{eq def Dc dan Cha Jha equiv bis}
    \Esp{ \sum_{t\in  \T}  q_t }=1,\;\;   1-\mu\le \chi_t \le
    1+\lambda
    \And\; S^2_t (1-\mu)E^{q}_t[{\bf 1}] \le E^{q}_t[\chi  S^2] \le E^{q}_t[{\bf 1}]
    S^2_t    (1+\lambda)
    \;,\;\;
     \end{equation}
where $E^{q}_t[\cdot]$ is defined as above with $ q$ in place of
$q^Z$. Recall from Theorem \ref{thm main} that there is some $Z
\in \Dc(K,\P)$ with values in $(0,\infty)^2$, so that
$\Qc(K,\P)\ne \emptyset$. \medskip

The set $\Qc(K,\P)$ coincides with the set of {\it
approximate martingale node-measure} defined in Chalasani and  Jha
\cite{CJ} (see Definition 6.3). More precisely, for $(\chi,q)$
$\in \Qc(K,\P)$,   the {\it node-measure} is defined by the map
    \b*
    A\in \Fc\x \Pc(\T) & \mapsto&
    \int \sum_{t\in \T} \1_{\{(\omega,t)\in A\}} d\P(\omega)
    \e*
where $\Pc(\T)$ is the collection of subsets of $\T$.  The term
$\chi$ does not appear in the formulation of the above paper
because the authors do no take into account the transaction costs
that are possibly paid when the option is exercised and the
hedging portfolio is liquidated.\medskip

Conversely, given $(\chi,q)$ $\in$  $\Qc(K,\P)$ it is clear that
we can find $Z \in   \Dc(K, \P )$ such that $(\chi^Z,q^Z)$ $=$
$(\chi,q)$. It follows from Corollary \ref{cor form duale} that
for $\vartheta$ $\in  L^0(\R^2 ;\F, \P)$ such that $\vartheta^1$
and $\vartheta^2$  are uniformly bounded from below, we have
    \be
    h(\vartheta;K)&:=&\inf\{x^1 \in \R~:~(x^1,0) \in \Gamma(\vartheta;K)\}\nonumber \\
    &=& \sup_{Z \in   \Dc(K,  \P )} E^{q^Z}_0[\vartheta^1+\chi^Z \vartheta^2S^2] \nonumber \\
    &=&\sup_{(\chi,q) \in   \Qc(K,  \P )} E^{q}_0[\vartheta^1+\chi  \vartheta^2S^2] \;.\label{eq formu duale node measures}
    \ee
This extends the dual formulation in terms of {\it node-measures}
obtained by Chalasani and  Jha \cite{CJ}, see Theorem 9.1, to the
general discrete time case where we take into account the
transaction costs that  are possibly paid when the option is
exercised. Since $\vartheta$ corresponds in our framework to a
claim labelled in units of the assets, the corresponding amounts
are given by the process $(\vartheta^1,\vartheta^2 S^2)$.

\subsection{``Randomized stopping times''}

In Theorem 9.1 of Chalasani and Jha \cite{CJ} one can also find an
equivalent formulation in terms of {\it randomized stopping
times}. A {\it randomized stopping times} $X$ is a non negative
$\F$-adapted process such that $\sum_{t\in \T}X_t=1$.  We denote
by $\Xc$ the set of all \textit{randomized stopping times}.
Observe that for a stopping time $\tau$, the process defined by
$X:=(\1_{\tau=t})_{t\in \T}$ belongs to $\Xc$.  Chalasani and Jha
\cite{CJ} show that there is a one to one correspondence between
node measures  and   pairs $(X,\Q)$ where $X$  is a {\it
randomized stopping times}   and  $\Q$ is a $\P$-absolutely
continuous probability measure (see Theorem 5.4 in \cite{CJ}).
This result can be easily extended to our framework as shown
below.\medskip

Given an adapted process $\chi$ such that $(1-\mu)\le \chi_t\le (1+\lambda)$  for
all $t\in \T$, the $\P$-equivalent measure $\Q$  is   called a $(\chi,X)$-{\it approximate
martingale measure} if for all $t \in \T$
    \be \label{eq contrainte sur randomized ta}
    S^2_t (1-\mu)X^{+}_t \le  \E^{\Q}\left[\sum_{t\le s\le T}  X_s \chi_s S^2_s ~|~\Fc_t\right]
     \le  X^{+}_t  S^2_t   (1+\lambda)
   \ee
where $X^{+}_t$ $:=$ $\E^{\Q}\left[\sum_{t\le s\le T}  X_s  ~|~\Fc_t\right]$ for $t \in \T$. Observe that $X^+_t=\sum_{t\le s\le T} X_s$ since $X$ is $\F$-adapted and   $\sum_{t\in \T}X_t=1$.
Denoting by $\Pc(X;K,\P)$ the associated set of pairs $(\chi,\Q)$ such that the above
inequalities hold, we then obtain as in Chalasani and  Jha
\cite{CJ} that
 \be\label{eq formulation random ta}
    h(\vartheta;K)&=& \sup_{X \in \Xc}\; \sup_{(\chi,\Q) \in \Pc(
    X;K,\P)} \E^{\Q}\left[\sum_{t \in \T} X_t
    \left(\vartheta^1_t+\chi_t \vartheta^2_t S^2_t\right)\right] =:h(\Xc)\;.
 \ee
Here again, the $\chi$ is added to the formulation of Chalasani
and  Jha \cite{CJ} to take into account the transaction costs that
are possibly paid when the option is exercised.\medskip

To obtain the last equality first, we argue in two steps :

1. First, take $(\chi,q)$ $\in$  $\Qc(K,\P)$ and define
 \b*
 N_t=\Esp{\sum_{t\le s\le T} q_s~|~\Fc_t}  &\And&  D_t=\Esp{\sum_{t\le s\le T} q_s~|~\Fc_{t-1}}
 \e*
for all $0\le t\le T$ with the convention $\Fc_{-1}=\Fc_{0}$.
Then, let the $\F$-adapted processes $H$ and $X$ be defined
inductively by $(H_0,X_0)=(1,q_0)$, $H_{t+1}=H_tN_{t+1} /
D_{t+1}\1_{D_{t+1}\ne 0}$ $+$ $H_t\1_{D_{t+1}= 0}$ and
$X_t=q_t/H_t$ for all $1\le t\le T$. Then, $H$ is a martingale
starting from $1$ and we can define the associated equivalent
probability measure $\Q $ by $d\Q/d\P=H_T$. Moreover, one easily
checks by using an inductive argument that $\sum_{j=0}^k
X_{T-j}=N_{T-k}/H_{T-k}$ for all $1\le k \le T$. Indeed, since
$q_T=N_T$ and $D_{T}= N_{T-1}-q_{T-1}$ one has
    \b*
    X_{T-1}+X_T&=&
    \frac{q_{T-1}}{H_{T-1}} + \frac{  N_{T-1}-q_{T-1}}{H_{T-1} }\1_{D_{T}\ne 0}
    + \frac{q_{T}}{H_{T-1} }\1_{D_{T}= 0}
    \= \frac{N_{T-1}}{H_{T-1}}
    \e*
since  on $\{D_T=0\}$ one has $q_T=0$ and therefore
$q_{T-1}=N_{T-1}$. Using the identities $D_{t}= N_{t-1}-q_{t-1}$,
$t\ge 1$, the same argument provides the above result.  For $k=T$,
this shows that $\sum_{t=0}^T X_{t}=N_{0}/H_{0}=1$ since $N_0$ $=$
$\Esp{\sum_{t\in \T} q_t}$ $=$ $1$. Hence, $X$ is a randomized
stopping time. Rewriting $q$ in terms of $(H,X)$ in \reff{eq def
Dc dan Cha Jha equiv bis}-\reff{eq formu duale node measures}, one
obtains \reff{eq contrainte sur randomized ta} and the expectation
entering in the definition   of $h(\Xc)$. This shows that
$h(\Xc)\ge h(\vartheta;K)$.

2. Conversely, observe that for  $X\in \Xc(K,\P)$   and  $(\chi,\Q)$ $\in$ $\Pc(X;K,\P)$,
then $(\chi,q)$ $\in \Qc(K,\P)$
with $q$ defined by $q_t:=X_tH_t/\Esp{\sum_{s\in \T} X_sH_s}$
and $H_t:=\Esp{d\Q/d\P~|~\Fc_t}$ for $t\in \T$.  In view of \reff{eq formu duale node measures},
this shows that $h(\Xc)\le h(\vartheta;K)$.

\section{Counter examples}

In this section, we first show that the duality relation
    \b*
    {\bf D}(K):~\Gamma(\vart;K)=\Theta(\vart;K)\;\pourtout \;
    \vart \in L^0(\R^d;\F)
    \e*
does not hold for a large class of $\Cc$-valued process $K\in \Kc$
(recall the definitions of $\Gamma$ and $\Theta$   in equations  \reff{GammaDefinition}  and
  \reff{eq mauvaise formualtion duale}). For $x$
$\in$ $\R^d$, let us define
    \be\label{eq def et for dual de c_t}
    c_t(x)&:=&\min\left\{c\in \R~:~c\ubf_1 - x \in K_t\right\}
    \;.
    \ee
In financial terms, $c_t(x)$ is the minimal amount, in terms of
the first asset, necessary to dominate $x$ in the sense of $K_t$
at time $t$. If the first asset is interpreted as a numeraire, it
corresponds to the {\it constitution value} of $x$ in terms of
this numeraire. Here, $\ubf_1$ stands for the $\R^d$ vector
$(1,0,\ldots,0)$.

\begin{Proposition}\label{prop marchepas} If there exists $x$ $\in$ $\R^d$ such that
    \begin{itemize}
    \item[(i)] $y-c_0(x)\ubf_1$ $\in$ $K^0_0$
$\Rightarrow$ $y-x$ $\in$ $K^0_0$ or $\Pro{y-x\in K_1}<1$
    \item[(ii)] $x-c_0(x)\ubf_1 \notin K_0$.
    \end{itemize}
Then, there exists $\vart$ such that
$\Theta(\vart;K)\ne\Gamma(\vart;K)$, i.e. ${\bf D}(K)$ is not
satisfied.
\end{Proposition}

The proof is postponed to the end of the section.

\begin{Remark}{\rm Condition (ii) means that there are directions
  with efficient frictions at time $0$. Condition (i) has the following
  interpretation. If a portfolio $y$ is   equivalent to the {\it
  constitution value} of $x$ then it dominates $x$ in the sense of $K_0$.
   However, since $x$ and $y$ have the same constitution value, $c_0(x)=c_0(y)$, it can
  not be {\it too large}. In particular, if it is not
  equivalent to $x$, then it can not dominate $x$ component by
  component.
  In that case, we assume that  there is randomness enough so that the
  probability that $y$ still dominates $x$ at time $1$ is less
  than $1$.
}
\end{Remark}

\begin{Remark}{\rm 1. If $K_0^0=\{0\}$ and $x \ne c_0(x)\ubf_1$ then (ii)
holds since by definition we already have $c_0(x)\ubf_1-x \in
K_0$. If we also assume that $\Pro{c_1(x)>c_0(x)}>0$ then (i) is
satisfied too.}
\end{Remark}

\begin{Example}{\rm
{\it 1. Efficient frictions:} consider the following cones
    \b*
    K_t &=& \left\{(x^1,x^2) \in \R^2~:
    x^1+(1+\lambda_t) x^2 \ge 0\;,\; x^1+(1-\mu_t)x^2\ge 0\right\}
    ,
    \e*
where $t\in \T:=\{0,1\}$, $\lambda_0 < \lambda_1$ and
$\mu_0,\mu_1$ $\in $ $(0,1)$. Observe that $K^0_0=\{0\}$. For $x =
(0,1)$, $c_0(x) = 1+\lambda_0< c_1(x)=1+\lambda_1$. Then, the
conditions of the remark above hold so that ${\bf D}(K)$ is not
true.

{\it 2. Partial frictions:} consider the preceding case where we
add an asset which has no transaction cost with the first one,
i.e.
    \b*
    K_t &=& \left\{(x^1,x^2,x^3) \in \R^3~:
    x^1+(1+\lambda_t) x^2 +x^3\ge 0\;,\; x^1+(1-\mu_t)x^2+x^3\ge
    0\right\}.
    \e*
We put $x=(0,1,0)$ so that assumption (ii) holds. We now check
(i). It is clear that if $y-c_0(x)\ubf_1$ $\in$ $K^0_0$ then $y-x$
$\notin$ $K^0_0$. Observe that $y=(y^1,0,y^3)$ with
$y^1+y^3=c_0(x)$, so $y^1 + (-1)(1+\lambda_1) + y^3<0$ which
implies that $y-x \notin K_1$. }
\end{Example}

On the contrary, we can also show that  ${\bf D}(K)$ does not
\textit{only hold} in the case where $K_t$ $ =$ $K^0_t+\R^d_+$, i.e. there is no transaction costs.

\begin{Proposition}\label{K0=0_dualite} There exists $(\Omega,\Fc,\P)$
and $K \in \Kc$ such that $NA(K)$ holds, $K^0_t=\{0\}$ for all
$t$, and such that for all $\vart \in L^0(\R^d;\F)$ we have
$\Theta(\vart;K)= \Gamma(\vart;K)$.
\end{Proposition}

\proof We take $\Omega$ trivial, i.e. $|\Omega|$ $=$ $1$ with
$\Fc_0=\Fc_T$ $=$ $\{\Omega,\emptyset\}$, and put $K=K_0$
constant. Then, $x$ $\in \Theta(\vart;K)$ reads $\sup\limits_{Z_t
\in K^*_t} Z_t\cdot (\vart_t-x)\le 0$, i.e. $x-\vart \in K_t$ for
all $t\in \T$.\ep \vspace{0.5cm}

This example shows that, for ${\bf D}(K)$ to be wrong, we need not
only to have non zero transaction costs but also enough randomness
in the direction where transaction costs are positive. \medskip

{\bf Proof of Proposition \ref{prop marchepas}: } Let $x$ be such
that $(i)-(ii)$ are satisfied. We consider the asset $\vart$
defined by $\vart_t= c_0(x) \ubf_1 \1_{\{t=0\}}+ x \1_{\{t>0\}}$.
>From the martingale property of $Z$,
    \b*
    \sup_{\tau \in \Tc(\T)}
    \Esp{Z_\tau \cdot \vart_\tau - Z_0 \cdot(c_0(x) \ubf_1)}
    &=&
    \sup_{\tau \in \Tc(\T)}
    \Esp{Z_\tau \cdot (x -c_0(x) \ubf_1)\1_{\{\tau>0\}}}
    \\&=&
    \Max\left\{0\;;\;Z_0 \cdot (x -c_0(x) \ubf_1)\right\}
    \e*
which is non positive by \reff{eq def et for dual de c_t}. Hence,
$c_0(x) \ubf_1 \in \Theta(\vart;K)$. If ${\bf D}(K)$ holds, then
there exists a portfolio $V \in$ $A(c_0(x) \ubf_1;K)$ such that
$V_0-c_0(x) \ubf_1\in K_0$ and therefore $V_0-c_0(x) \ubf_1\in
K_0^0$.  By $(i)$ there is two cases. If $V_0 - x\in K^0_0$, then
$x -c_0(x) \ubf_1$ $\in$ $K^0_0 \subset K_0$ which is a
contradiction of $(ii)$. If $\Pro{V_0-x \in K_1}<1$,  we can not
have $V_1-x=V_0+\xi_1-x \in K_1$  with $\xi_1 \in -K_1$. \ep

\section{Proofs}

\no In this section, we first provide the proof of Theorem
\ref{thm main}. It follows from standard arguments based on the
Hahn-Banach separation theorem. For ease of notations, we simply
write $A(K)$ and $A^s(K)$ in place of $A(0;K)$ and $A^s(0;K)$. We
denote by $\Lc^0$ the set of $\F$-adapted processes with values in
$\R^d$ and by $\Lc^1(\tilde \P)$ (resp. $\Lc^\infty$) the subset
of these elements which are $\tilde \P$-integrable, $\tilde \P
\sim \P$, (resp. bounded). Observe that $\Lc^0$ (resp.
$\Lc^\infty$) can be identified as a subset of $L^0(\R^{d \times
(T+1)};\Fc)$ (resp. $L^\infty(\R^{d \times (T+1)};\Fc)$, the set
of bounded elements of $L^0(\R^{d \times (T+1)};\Fc)$).

\begin{Proposition}\label{prop separateur} Let the conditions of  Theorem \ref{thm
main} hold.  Then, for all $\tilde \P\sim \P$, there is some $Z
\in \Dc(K;\tilde \P )\cap \Lc^\infty$ with values in
$(0,\infty)^d$  such that
    \b*
    \sup_{\vartheta \in A^s(K)\cap \Lc^1(\tilde \P)}
    \Espt{
    \sum_{t \in \T} Z_t \cdot \vartheta_t
    }
    &\le&
    0\;.
    \e*
\end{Proposition}

\proof Since $A^s(K) \cap \Lc^1(\tilde \P)$ is closed in  $
L^1(\R^{d\times (T+1)};\Fc,\tilde \P)$ (when identified with a
subset of $L^1(\R^{d\times (T+1)};\Fc,\tilde \P)$) and convex, it
follows from the Hahn-Banach separation theorem, $NA(K)$ and the
fact that $A^s(K)\cap \Lc^1(\tilde \P)$ is a cone which contains
$-\Lc^\infty$, that, for all $\phi \in L^1(\R^d_+;\F,\tilde
\P)\setminus \{0\}$,  there is some $\eta^\phi=(\eta^\phi_t )_{t
\in \T} \in L^\infty(\R_+^{d\times (T+1)};\Fc)$ such that
    \b*
    \sup_{\vartheta \in A^s(K)\cap \Lc^1(\tilde \P)}
    \Espt{
    \sum_{t \in \T} \eta^\phi_t \cdot \vartheta_t
    }
    &\le& 0 <
    \Espt{
    \sum_{t \in \T} \eta^\phi_t \cdot \phi_t
    }\;.
    \e*
By  possibly replacing $\eta^\phi_t$ by
$\Esp{\eta^\phi_t~|~\Fc_t}$, we can assume that $\eta^\phi$ is
$\F$-adapted. Then, using a standard exhaustion argument, one can
find some $\eta \in \Lc^\infty$ with values in $(0,\infty)^{d}$
such that
    \be\label{eq espt nega}
    \sup_{\vartheta \in A^s(K)\cap \Lc^1(\tilde \P)}
    \Espt{
    \sum_{t \in \T} \eta_t \cdot \vartheta_t
    }
    &\le& 0 \;.
    \ee
Fix  some arbitrary $\xi \in \Ac(K)\cap \Lc^\infty$, so that
$V^{0,\xi}$ $\in$ $A^s(K) \cap \Lc^1(\tilde \P)$. Since
    \b*
    \sum_{t \in \T} \eta_t \cdot V^{0,\xi}_t
    &=&
    \sum_{t \in \T}   \xi_t \cdot \left(\sum_{s=t}^T  \eta_s\right)
    \e*
we deduce from the above inequality that
    \b*
    \sup_{\xi \in \Ac(K)\cap \Lc^\infty}
    \Espt{
    \sum_{t \in \T} \bar \eta_t \cdot \xi_t
    }
    &\le& 0\;,
    \e*
where we defined
    \b*
    \bar \eta_t&:=& \Espt{\sum_{s=t}^T \eta_s~|~\Fc_t}  \; \;\;t\in \T\;.
    \e*
This shows that $\bar \eta_t \in K^*_t$ for all $t\in \T$.
For an arbitrary bounded  element  $\xi_t$ in $L^0(K_t;\Fc_t)$,
the process $V^{0,\xi}_s = -\1_{\{s=t\}}\xi_t$, $s \in \T$,
belongs to $A^s(K)$. In view of \reff{eq espt nega}, this implies
that $\eta_t \in K^*_t$.
 \ep

\begin{Proposition}\label{prop cond suff est dans A} Let the conditions of  Theorem \ref{thm
main} hold. Fix $\tilde \P\sim \P$ and $\vart\in \Lc^1(\tilde
\P)$. If
    \b*
    \Espt{\sum_{t=0}^T \vart_t \cdot Z_t} &\le&
    0\;
    \e*
for all  $Z \in \Dc(K,\tilde \P)$ such that $\vart \cdot Z$ $\in$
$\Lc^1(\tilde \P)$, then   $\vart \in A^s(K)$.
\end{Proposition}

\proof Since $A^s(K)\cap \Lc^1(\tilde \P)$ is closed and convex,
if $\vart$ $\notin A^s(K)$, we can find some $\eta=(\eta_t)_{t \in
\T} \in L^\infty(\R^{d\times (T+1)};\Fc)$ such that
    \b*
    \sup_{\tilde \vart \in A^s(K)\cap \Lc^1(\tilde \P)}
    \Espt{\sum_{t=0}^T \tilde \vart_t \cdot \eta_t} &<& \Espt{\sum_{t=0}^T \vart_t \cdot \eta_t}
    \;.
    \e*
By the   same arguments as in the proof of Proposition \ref{prop
separateur}, we can assume that  $\eta$ is $\F$-adapted and show
that $\eta_t \in K^*_t$ and $\bar \eta_t\in K^*_t$  for all
$t\in \T$. Hence, $\eta \in \Dc(K,\tilde \P)$ which   leads to a contradiction.
 \ep \vspace{0.5cm}

{\bf Proof of Theorem \ref{thm main}} 1. In view Proposition
\ref{prop cond suff est dans A}, the implication (ii)
$\Rightarrow$ (i) is obtained by considering $\tilde \P$ with
density with respect to $\P$ defined by $H/\Esp{H}$ with
$H:=\exp(-\sum_{t\in \T} \|\vart_t\|)$.

\no 2. It is clear that (ii) implies (iii) while the reverse
implication follows from the fact that $Z \in \Dc(K,  \P)$ if and
only if $\tilde H Z \in \Dc(K,  \tilde \P)$ where $\tilde
H_t:=\Esp{d \P/d\tilde\P~|~\Fc_t}$.

\no 3. The last implication  (i) $\Rightarrow$ (ii) is trivial.
Indeed, recall that, for $\xi \in \Ac(K)$,
    \b*
    \Esp{\sum_{t \in \T} Z_t \cdot V^{0,\xi}_t}
    &=&
    \Esp{\sum_{t \in \T} \bar Z_t \cdot \xi_t}\;.
    \e*
Since $\bar Z_t \in L^0(K^*_t;\Fc_t)$ and $\xi_t \in
L^0(-K_t;\Fc_t)$, the last term is non-positive. Moreover,
$V^{0,\xi}_t - \vart_t$ $\in$ $L^0(K_t;\Fc_t)$ implies $ Z_t \cdot
V^{0,\xi}_t$ $\ge$ $ Z_t\cdot \vart_t$.
 \ep \vspace{0.5cm}

\no We now provide the proof of Proposition \ref{lem closure}. The
following Lemma can be found in  Kabanov  and Stricker
\cite{teachers note}.

\begin{Lemma}\label{lem liminf finie} Set $\Gc \subset \Fc$ and $E$ be a closed subset of $\R^{d}$.
Let $(\eta^n)_{n\ge 1}$ be a sequence in $L^0(E;\Gc)$. Set $\tilde
\Omega$ $:=$ $\{\liminf_{n\to \infty} \|\eta^n\| < \infty\}$.
Then, there is an increasing sequence of random variables
$(\tau(n))_{n\ge 1}$ in $L^0(\N;\Gc)$ such that $\tau(n)\to
\infty$  and , for each $\omega \in \tilde \Omega$,
$\eta^{\tau(n)}(\omega)$ converges to some $\eta^*(\omega)$ with
$\eta^* \in L^0(E;\Gc)$.
\end{Lemma}

\no {\bf Proof of Proposition \ref{lem closure}. } We use an
inductive argument. For $t\in \T$, we denote by $\Sigma_t$ the set
of processes $\vartheta$ $\in \Lc^0$ such that
    \b*
    \exists\; \xi \in \Ac(K) \;\mbox{ s.t.}\;\; \sum_{s=t}^\tau \xi_s - \vartheta_\tau &\in & K_\tau
    \;\;\pourtout\;t\le \tau \le T\;.
    \e*
Clearly, $\Sigma_T$ is closed in measure. Assume that
$\Sigma_{t+1}$ is closed and let $\vartheta^n$ be a sequence in
$\Sigma_{t}$ such that $\vartheta^n_s\to \vartheta_s$  for
$t\le s\le T$. Let $\xi^n$ $\in \Ac(K)$ be such that
    \b*
    \sum_{s=t}^\tau \xi^n_s - \vartheta^n_\tau &\in & K_\tau
    \;\;\pourtout\;t\le \tau \le T\;.
    \e*
Set $\tilde \Omega$ $:=$ $\left\{ \liminf_{n\to \infty}
\|\xi^n_t\|<\infty\right\}$. Since $\tilde \Omega \in \Fc_t$, we
can work separately on $\tilde \Omega$ and $\tilde \Omega^c$.

\no 1. If $\Pro{\tilde \Omega}=1$, after possibly passing to a
random sequence (see Lemma \ref{lem liminf finie}), we can assume
that $\xi^n_t$ converges $\Pas$ to some $\xi_t \in
L^0(-K_t;\Fc_t)$. Since
    \b*
    \sum_{s=t+1}^\tau \xi^n_s - (\vartheta^n_\tau-\xi^n_t) &\in & K_\tau
   \;\;\pourtout\;t+1\le \tau \le T\;,
    \e*
and $\Sigma_{t+1}$ is closed, we can find some $\tilde \xi$ $\in$
$\Ac(K)$ such that
    \b*
    \sum_{s=t+1}^\tau \tilde \xi_s - (\vartheta_\tau-\xi_t) &\in & K_\tau
    \;\;\pourtout\;t+1\le \tau \le T\;.
    \e*
This shows that $\vartheta \in \Sigma_{t}$.

\no 2. If $\Pro{\tilde \Omega}<1$, then we can assume without loss
of generality that $\Pro{\tilde \Omega}=0$. Following line by line
the proof of Lemma 2 in Kabanov et al. \cite{KSR} and using the
$K_s$'s closure property, we can find some $\hat \xi \in \Ac(K)$
with $\| \hat \xi_t \| =1$ such that
    \b*
    \kappa_\tau:=\sum_{s=t}^\tau \hat \xi_s  &\in & K_\tau\;\;\;\pourtout\;t\le \tau \le T\;.
    \e*
Since that $0 = \sum_{s=t}^\tau \hat\xi_s - \kappa_\tau =
\sum_{s=t}^{\tau-1} \hat\xi_s+(\hat\xi_\tau - \kappa_\tau)$ and
that $\hat\xi_\tau$ and $-\kappa_\tau$ $\in -K_\tau$,
we deduce by \reff{eq_ha} that $\hat\xi_\tau - \kappa_\tau \in
K_\tau^0$. Therefore,
    \be\label{eq somme des hat xi dans K0}
   \hat \xi_\tau \in K^0_\tau \;\And\; \kappa_\tau=\sum_{s=t}^\tau \hat \xi_s  &\in & K^0_\tau
   \;\;\pourtout\;t\le \tau \le T\;.
    \ee
Since $\|\hat \xi_t\|=1$, there is a partition of $\tilde \Omega$
into disjoint subsets $\Gamma_i \in \Fc_t$ such that $\Gamma_i
\subset \{(\hat \xi_t)^i\ne 0\}$ for $i=1,\ldots,d$. We then
define
    \b*
     \check \xi^n_s&=& \sum_{i=1}^d \left(\xi^n_s - \beta^{n,i}_t \hat \xi_s\right)\1_{\Gamma_i} \;\;\;s\in \T\;
    \e*
with $\beta^{n,i}_t$ $=$ $(\xi^{n}_t)^i/(\hat \xi_t)^i$  on
$\Gamma_i$,   $i=1,\ldots,d$. In view of \reff{eq somme des hat xi
dans K0} and definition of $\xi^n$, we have
    \b*
    \sum_{s=t}^\tau \check \xi^n_s -  \vartheta^n_\tau  &\in & K_\tau
    \;\;\pourtout\;t\le \tau \le T\;,
    \e*
since $K_\tau-K^0_\tau\subset K_\tau, \;\;\tau\in \T$. We can then
proceed as in Kabanov et al. \cite{KSR} and obtain the required
result by repeating the above argument with $(\check \xi^n)_{n\ge
1}$ instead of $(\xi^n)_{n\ge 1}$ and by iterating this procedure
a finite number of times.
 \ep

\end{document}